\documentclass[final,1p, times]{elsarticle}

\usepackage{amssymb,amsmath}

\usepackage{lineno}

\newtheorem{thm}{Theorem}[section]

\newtheorem{lemma}[thm]{Lemma}

\newproof{pf}{Proof}
\newcommand{\bF}{\mathbb{F}_{q}}

\begin{document}

\begin{frontmatter}



\title{Suborbits of a  point stabilizer in the orthogonal group on the
last subconstituent of orthogonal dual polar graphs }

\author[FL]{Fenggao Li}\ead{fenggaol@163.com}
\author[KW]{Kaishun Wang\corref{cor}}
\ead{wangks@bnu.edu.cn}
\author[JG]{Jun Guo}\ead{guojun-lf@163.com}
\author[JM]{Jianmin Ma}
\ead{Jianmin.Ma@emory.edu}
\cortext[cor]{Corresponding author}

\address[FL]{College of Mathematics, Hunan Institute of Science and Technology, Yueyang  414006,   China}
\address[KW]{Sch. Math. Sci. \& Lab. Math. Com. Sys.,
Beijing Normal University, Beijing  100875, China}
\address[JG]{Math. and Inf. College, Langfang Teachers' College, Langfang  065000,  China }
\address[JM]{Oxford College of Emory University, Oxford, GA 30054, USA}

\begin{abstract}
As one of the serial papers on suborbits of point stabilizers in classical groups on the
last subconstituent of dual polar graphs, the corresponding problem
for orthogonal dual polar graphs over a finite field of odd characteristic
is discussed in this paper. We determine all the suborbits of a point-stabilizer in the orthogonal group on the last subconstituent, and calculate the length of each suborbit. Moreover,
we discuss the quasi-strongly regular graphs and the association schemes based on the last subconstituent, respectively.
\end{abstract}

\begin{keyword}
orthogonal group \sep suborbit \sep dual polar graph \sep
subconstituent \sep   quasi-strongly regular graph \sep association scheme


\MSC[2010] 05B125 \sep 05C25

\end{keyword}

\end{frontmatter}

\section{Introduction}

Let $\bF$ be a finite field with $q$ elements, where $q$ is an odd prime power.
Let $\bF^{n}$ be the row vector space of dimension $n$ over $\bF$.
The set of all $m\times n$ matrices over $\bF$ is denoted by $M_{mn}(\bF)$, and $M_{nn}(\bF)$ is denoted by $M_n(\bF)$ for simplicity.
For any matrix $A=(a_{ij})\in M_{mn}(\bF)$, we denote the transpose of $A$ by $A^t$.

Let $n=2\nu+\delta$, where $\nu$ is a non-negative integer and $\delta=0,1$ or 2. Suppose
$$
\begin{array}{rl}
S_{2\nu+\delta,\Delta}\!=\!\left(\!\begin{array}{ccc}\!0\!&I^{(\nu)}\!&\!\\
                                      \!I^{(\nu)}\!&0\!&\!\\
                                      \!&\!&\!\Delta\!\end{array}\!\right),
&
\Delta=\left\{\begin{array}{ll}\phi\,\mbox{(disappear)},&\mbox{if }\delta=0,\\
                  (1)\mbox{ or } (z),      &\mbox{if } \delta=1,\\
                  \left(\!\begin{array}{cc}\!1\!&\!\\[-0.05cm]
                        \!\!&\!-z\!\end{array}\!\right),&\mbox{if } \delta=2,
\end{array}\right.
\end{array}
$$
where $z$ is a fixed non-square element of $\bF$ such that $1-z$ is a non-square element.
When $\delta=1$ or 2, $\Delta$ is definite in the sense that for
any row vector $x\in\bF^{\delta}$, $x{\Delta}x^T=0$ implies $x=0$. Note that the set
$$
\left\{\,T\in GL_{2\nu+\delta}(\bF)\,|\,TS_{2\nu+\delta,\Delta}T^t=S_{2\nu+\delta,\Delta}\,\right\}
$$
forms a subgroup of $GL_{2\nu+\delta}(\bF)$, called the {\em orthogonal group} of degree
$n=2\nu+\delta$ with respect to $S_{2\nu+\delta,\Delta}$ over $\bF$, denoted by
$O_{2\nu+\delta,\Delta}(\bF)$.
The group $O_{2\nu+\delta,\Delta}(\bF)$ acts on $\bF^{2\nu+\delta}$ by the matrix multiplication.
$\bF^{2\nu+\delta}$ together with this action  is called
the $(2\nu+\delta)$-dimensional {\em orthogonal space} over $\bF$ with respect to
$S_{2\nu+\delta,\Delta}$.
A matrix representation of a subspace $P$ is a matrix whose rows form a basis for $P$.
When there is no danger of confusion, we use the same symbol to denote
a subspace  and its matrix representation.
An $m$-dimensional subspace $P$ of $\bF^{2\nu+\delta}$ is called {\em totally isotropic}
if $PS_{2\nu+\delta,\Delta}P^t=0$. It is well-known that
maximal totally isotropic subspaces of $\bF^{2\nu+\delta}$ are of dimension $\nu$.

Let $G$ be a  group acting transitively on
a finite set $X$. For a fixed element $a\in X$, the
stabilizer $G_a$  is not transitive on $X$ in general. The
orbits of $G_a$ on $X$ are said to be {\em suborbits}, and the number of
such suborbits is the {\em rank} of this action.
H. Wei and Y. Wang \cite{ww,Wei1,wei2} studied the suborbits of
the transitive set of all totally isotropic subspaces  under finite
classical groups. We discussed these
problems in singular classical spaces  in \cite{gwl, wgl}.

Dual polar graphs are famous distance-regular graphs and have been well studied (\cite{BI,BCN,Wan0}).
The {\em orthogonal dual polar graph} $\Gamma$ (on the orthogonal space $\bF^{2\nu+\delta}$) has as vertices the maximal totally
isotropic subspaces;  two vertices $P$ and $Q$ are adjacent if and only
if $\dim (P\cap Q) =\nu-1$. It is well-known that
$\Gamma$ is of diameter $\nu$.
For any vertex $P$ of $\Gamma$, the $i$th subconstituent
$\Gamma_i(P)$ with respect to $P$ is the induced graph on
the set of vertices at distance $i$ from $P$ in $\Gamma$.
A. Munemasa \cite{M} initiated the study of the subconstituents of dual polar
graphs in the orthogonal spaces, and characterized the first and last subconstituents.
Subsequently, Y. Wang, F. Li and Y. Huo \cite{Li1, Li2, Wang2, Wang1} characterized all
the subconstituents of dual polar graphs under finite classical groups,
and proved that for any vertex
$P$ of the dual polar graph $\Gamma$ in the $(2\nu+\delta)$-dimensional classical space (where $\delta=0,1\,\,\mathrm{or}\,\,2$),
the $m$th subconstituent $\Gamma_{m}(P)$ is isomorphic to
${\nu\brack m}_{q}\cdot{\cal G}^{(m,\delta)}$, where ${\cal G}^{(m,\delta)}$
is the graph with the vertex set consisting of the matrices $(X\ \ Z)$ such that
$$
\left\{\begin{array}{lll}
X+\overline{X}^t+ZI^{(\delta)} \overline{Z}^t=0& \mbox{the unitary case},\\
X+X^t+Z\Delta Z^t=0& \mbox{the orthogonal case of odd characteristic},
\end{array}\right.
$$
where $X\in M_m(\bF),\,Z\in M_{m\delta}(\bF)$; and two vertices $(X\ \ Z)$ and $(X_1\ \ Z_1)$ are adjacent if and only if
$(X - X_1\,\,Z - Z_1) $ is of rank $1$.
Note that the mapping
$$
(X\ \ Z)\mapsto (X\ \ I^{(m)}\ \ Z)
$$
is an isomorphism from ${\cal G}^{(m,\delta)}$ to the last subconstituent
of the corresponding dual polar graph in the classical space
$\bF^{2m+\delta}$. Therefore, the
study of subconstituents of a dual polar graph may be reduced to that of the
last subconstituent. In \cite{lwgm} we studied
the suborbits of a  point-stabilizer in the unitary group on the
last subconstituent of Hermitean dual polar graphs. In this paper
we discuss the corresponding problem
for orthogonal dual polar graphs over a finite field of odd characteristic.

Let $\Gamma$ be the orthogonal dual polar graph. It is
well-known that a point-stabilizer of $P$ of $\Gamma$ in $O_{2\nu+\delta,\Delta}(\bF)$
is transitive on the last subconstituent of
$\Gamma$. In Section 2 we determine all the suborbits of this action, and
calculate the rank and the lengths of these suborbits.
As two applications of our results, in Sections 3 and 4, we discuss the quasi-strongly regular graphs and the association schemes based on the last subconstitute of $\Gamma$, respectively.

\section{Suborbits}

Let $\Gamma$ be the dual polar graph in the orthogonal space $\bF^{2\nu+\delta}$.
Note that the last subconstituent $\Lambda$ is a coclique when $\delta=0$ (see \cite{Li1}), and
$\Lambda$ is studied in \cite{Li2} when  $\delta=1$. So the case $\delta =2$ is the
main objective of this paper.

Denote by $[X_1,X_2,\dots,X_t]$ the
block diagonal matrix whose blocks along the main diagonal are
matrices $X_1,X_2,\dots,X_t$, by $\mathcal{A}_{2r}=[\mathcal{A}_2,\ldots,\mathcal{A}_2]$ the $2r\times 2r$ matrix of rank $2r$ in which
 $\mathcal{A}_2=\left(\!\begin{array}{cc}
             0&1\\
             -1&0\\
            \end{array}
      \!\right)$,
 and by $(A\; B \dots C)$ the matrix whose
block entries are  $A, B, \dots, C$. Suppose $I^{(m)}$ denotes the identity matrix
of order $m$, and $0^{(p,q)}$ denotes the zero matrix of order $p$ by $q$ or
$0^{(p)}$ when $p=q$.

We now study  the suborbits of the stabilizer of each vertex $P_0$ in $O_{2\nu+2,\Delta}(\bF)$ on $\Lambda$.
Since $O_{2\nu+2,\Delta}(\bF)$ acts transitively on the subspaces of the same type,
we may choose $P_0=(I^{(\nu)}\ \ 0^{(\nu)}\ \ 0^{(\nu,2)})$.
By \cite{Li1}, $\Lambda$ consists of subspaces of form $(A\ \ I^{(\nu)}\ \ Z)$,
where $A\in M_{\nu}(\bF)$ and $Z\in M_{\nu 2}(\bF)$ satisfy $A+A^t+Z\Delta Z^t=0$.
Let $G_0$ be the stabilizer of $P_0$ in $O_{2\nu+2,\Delta}(\bF)$.
Then  $G_0$ consists of matrices
of  the following form:
$$
\left(\begin{array}{ccc}
                T_{11}&0&0\\
                T_{21}&(T_{11}^{t})^{-1}&T_{23}\\
                -S\Delta T_{23}^{t}T_{11}&0&S\end{array}\right),
$$
where $T_{11}\in GL_\nu(\mathbb{F}_q)$, $T_{21}\in M_\nu(\mathbb{F}_q)$, $T_{23}\in
M_{\nu 2}(\mathbb{F}_q)$, $S\in O_{2\times 0+2,\Delta}(\mathbb{F}_q)$ and
$$
(T_{11}^{t})^{-1}T_{21}^{t}+T_{21}T_{11}^{-1}+T_{23}\Delta T_{23}^{t}=0.
$$
It is well-known that $G_0$ acts transitively on $\Lambda$.
For any $P_1\in \Lambda$, the suborbits of $G_0$  are just the orbits of the
point-stabilizer of $P_1$ in $G_0$ on $\Lambda$.
Let  $P_{1}=(0^{(\nu)}\ \ I^{(\nu)}\ \ 0^{(\nu,2)})\in \Lambda$ and
$G_{01}$ be the stabilizer of  $P_0$ and $P_1$ in
$O_{2\nu+2,\Delta}(\bF)$. Then $G_{01}$ consists of matrices
of  the following form:
\begin{equation}\label{Eq1}
[T,\,(T^{t})^{-1},\,S],
\end{equation}
where $T\in GL_{\nu}(\bF)$ and $S\in O_{2\times 0+2,\Delta}(\bF)$.
The action of $O_{2\nu+2,\Delta}(\bF)$ on $\bF^{2\nu+2}$
induces an action  $G_{01}$ on $\Lambda$:
$$
\begin{array}{ccc}
 \Lambda\times G_{01} & \longrightarrow & \Lambda\\
                         ((A\ \ I^{(\nu)}\ \ Z),[T,\,(T^t)^{-1},\,S]) & \longmapsto
                         &(T^tAT\ \ I^{(\nu)}\ \ T^tZS).
\end{array}
$$

Denote by ${\sf K}_n$ the set of all $n\times n$ alternate matrices over $\mathbb{F}_q$.
In order to determine the orbits of $G_{01}$ on $\Lambda$, we need to introduce
an action on ${\sf K}_{\nu}$.
For $i=1,2$, let $\mathcal{O}_{i}$ denote the set of all matrices of the form
\begin{equation}\label{Eq2}
\left(\begin{array}{cc}
T_{11}&0\\
T_{21}&T_{22}
\end{array}\right),
\end{equation}
where $T_{11}\in
GL_{i}(\bF)$, $T_{21}\in M_{\nu-i,i}(\bF)$ and $T_{22}\in
GL_{\nu-i}(\bF)$. Then
$\mathcal{O}_{i}$ is a subgroup of $GL_{\nu}(\bF)$, and
there is an
action of $\mathcal{O}_{i}$ on ${\sf K}_{\nu}$:
$$
\begin{array}{ccc}
  {\sf K}_{\nu}\times \mathcal{O}_{i} & \longrightarrow & {\sf K}_{\nu}\\
                         (X,\,T) & \longmapsto     & T^{t}XT.
\end{array}
$$
Note that $\{0^{(\nu)}\}$ is the trivial orbit of $\mathcal{O}_{i}$ on ${\sf K}_{\nu}$ for $i=1,2$.

\begin{thm}\label{thm2.1}
{\rm (i)} The nontrivial orbits of $\mathcal{O}_{1}$ on ${\sf K}_{\nu}$ have the following representatives:
\begin{equation}\label{Eq3}
[0,\mathcal{A}_{2r},0^{(\nu-2r-1)}]\,\,(1\leq r\leq \lfloor(\nu-1)/2\rfloor),\quad [\mathcal{A}_{2r},0^{(\nu-2r)}]\,\,(1\leq r\leq \lfloor\nu/2\rfloor).
\end{equation}

{\rm (ii)} The nontrivial orbits of $\mathcal{O}_{2}$ on ${\sf K}_{\nu}$ have the following representatives:
\begin{equation}\label{Eq4}
\begin{array}{lll}
&[0,\mathcal{A}_{2r},0^{(\nu-2r-1)}]\,(1\leq r\leq \lfloor(\nu-1)/2\rfloor),&
[0^{(2)},\mathcal{A}_{2r},0^{(\nu-2r-2)}]\,(1\leq r\leq \lfloor(\nu-2)/2\rfloor),\\[0.1cm]
&[\mathcal{A}_{2r},0^{(\nu-2r)}]\,(1\leq r\leq \lfloor\nu/2\rfloor),&
[\mathcal{K},\mathcal{A}_{2r-4},0^{(\nu-2r)}]\,(2\leq r\leq \lfloor\nu/2\rfloor),\\[0.1cm]
\end{array}
\end{equation}
where
$\mathcal{K}=\left(\begin{array}{cc}
 0&I^{(2)}\\
 -I^{(2)}&0\\
\end{array}\right)$.
\end{thm}
\begin{pf}
We only prove (ii), and (i) can be treated similarly.
Let $X\in {\sf K}_{\nu}$ with rank $2r>0$. Write
$$
X=\left(\!\begin{array}{cc}
        x\mathcal{A}_2&X_{12}\\
        -X_{12}^{t}&X_{22}\\
\end{array}\!\right),
$$
where $X_{12}\in M_{2,\nu-2}(\mathbb{F}_q)$ and $X_{22}\in{\sf K}_{\nu-2}$.
Then
$\mathrm{rank}\,X_{22}=2(r-i)$, $i=0,1$ or 2. Hence there is a $T_{11}\in
GL_{\nu-2}(\bF)$ such that
$T_{11}^{t}X_{22}T_{11}=[\mathcal{A}_{2(r-i)},0^{(\nu-2r+2i-2)}]$. Let
$$
X_{12}T_{11}=\bordermatrix{ &_{2r-2i}&_{\nu-2r+2i-2}\cr
 &Y_{12} & Y_{13}}\,\,\,\,\mbox{and}\,\,\,
T=\left(\begin{array}{cc}
   I^{(2)}&\\
   &T_{11}
\end{array}\right)\!
\left(\begin{array}{ccc}
   I^{(2)}&&\\
   -\mathcal{A}_{2(r-i)}Y_{12}^{t}&I^{(2r-2i)}&\\
   && I^{(\nu-2r+2i-2)}\\
\end{array}\right).
$$
Then $T\in\mathcal{O}_2$ and
$$
T^{t}XT=\left(\begin{array}{ccc}
              x_1\mathcal{A}_2&0&Y_{13}\\
              0&\mathcal{A}_{2(r-i)}&0\\
              -Y_{13}^{t}&0&0^{(\nu-2r+2i-2)}
\end{array}\right),
$$
where $x_1\mathcal{A}_2=x\mathcal{A}_2-Y_{12}\mathcal{A}_{2(r-i)}Y_{12}^{t}$.

{\em Case 1:} $\mathrm{rank}\,X_{22}=2r$.  Then $x_1=0$, $Y_{13}=0$ and
$T^{t}XT=[0^{(2)},\mathcal{A}_{2r},0^{(\nu-2r-2)}]$.

{\em Case 2:} $\mathrm{rank}\,X_{22}=2(r-1)$. Then
$\mathrm{rank}\,Y_{13}=0$ or 1. If $\mathrm{rank}\,Y_{13}=0$, then $x_1\ne 0$, $T_{1}=[x_1^{-1},I^{(\nu-1)}]\in\mathcal{O}_2$ and
$T_1^{t}T^{t}XTT_1=[\mathcal{A}_{2r},0^{(\nu-2r)}]$.
If $\mathrm{rank}\,Y_{13}=1$,
then there exists a $T_{12}\in GL_{2}(\mathbb{F}_q)$ and $T_{13}\in GL_{\nu-2r}(\mathbb{F}_q)$ such that
$$
T_{12}Y_{13}T_{13}=\left(\begin{array}{cc}0&0^{(\nu-2r-1)}\\ 1&0^{(\nu-2r-1)}\end{array}\right).
$$
Let $x_2\mathcal{A}_{2}=x_1T_{12}\mathcal{A}_2T_{12}^t$ and
$$
T_{2}=\left(\begin{array}{ccc}
       T_{12}^t&&\\
       &I^{(2r-2)}&\\
       &&T_{13}
\end{array}\right)\!
\left(\begin{array}{cccc}
       I^{(2)}&&&\\
       &0&I^{(2r-2)}&\\
       &1&0&\\
       &&&I^{(\nu-2r-1)}\\
 \end{array}\right)\left(\begin{array}{cccc}
       1&&&\\
       &1&&\\
       x_{2}&&1&\\
       &&&I^{(\nu-3)}\\
 \end{array}\right).
$$
Then $TT_{2}\in\mathcal{O}_2$ and
$(TT_{2})^{t}X(TT_{2})=[0,\mathcal{A}_{2r},0^{(\nu-2r-1)}]$.

{\em Case 3:} $\mathrm{rank}\,X_{22}=2(r-2)$. Then
$\mathrm{rank}\,Y_{13}=2$; and so there exists a $T_{14}\in GL_{\nu-2r+2}(\mathbb{F}_q)$ such that
$Y_{13}T_{14}=(I^{(2)}\ \ 0^{(2,\nu-2r)})$.
Let
$$
T_{3}=\left(\begin{array}{cc}
       I^{(2r-2)}&\\
       &T_{14}
 \end{array}\right)\!
\left(\begin{array}{cccc}
       I^{(2)}&&&\\
       &0&I^{(2r-4)}&\\
       &I^{(2)}&0&\\
       &&&I^{(\nu-2r)}\\
 \end{array}\right)
 \left(\begin{array}{cccc}
       1&&&\\
       &1&&\\
       &-x_1&1&\\
       &&&I^{(\nu-3)}\\
 \end{array}\right).
$$
Then $TT_{3}\in\mathcal{O}_2$ and
$(TT_{3})^{t}X(TT_{3})=[\mathcal{K},\mathcal{A}_{2(r-2)},0^{(\nu-2r)}]$.

Note that matrices of difference ranks can not be in the same orbit.
Now we show that any two distinct matrices in (\ref{Eq4}) cannot fall into the same orbit of $\mathcal{O}_{2}$. Otherwise,
there exists a $T\in\mathcal{O}_{2}$, which is of the form (\ref{Eq2}), carrying
$[0,\mathcal{A}_{2r},0^{(\nu-2r-1)}]$ to $[0^{(2)},\mathcal{A}_{2r},0^{(\nu-2r-2)}]$, then
$T_{22}^t[0,\mathcal{A}_{2(r-1)},0^{(\nu-2r-1)}]T_{22}=[\mathcal{A}_{2r},0^{(\nu-2r-2)}]$, which is impossible since $T_{22}$ is nonsingular.
Similarly, the left cases may be handled.

By above discussion, the desired result follows. \hfill $\Box$
\end{pf}

To determine the orbits of $G_{01}$ on $\Lambda$,  we need the following two lemmas.

\begin{lemma}~\label{lem2.2}
Let $a,b\in\bF$ with $(a,b)\ne (0,0)$. Then there exists a
$T\in O_{2\times 0+2,\Delta}(\bF)$ such that the subspace $(a,b)T$ has the matrix representation of the form
$(1,0)$ or $(1,1)$ corresponding to $a^2-zb^2$ is a square element or not, respectively.
\end{lemma}
\begin{pf}
Note that $(a,b)$ is of type $(1,1,0,1)$ or $(1,1,0,z)$ in $\bF^{2\times 0+2}$ corresponding to $a^2-zb^2$ being a square element or not, respectively. The result follows from \cite[Theorem~6.4]{Wan}. \hfill $\Box$
\end{pf}

\begin{lemma}~\label{lem2.3}
Any element of $O_{2\times 0+2,\Delta}(\bF)$ has one of the following forms
\begin{equation}\label{Eq5}
\left(
    \begin{array}{cc}
      x&y\\
      yz&x\\
    \end{array}
  \right),\,\,
\left(
    \begin{array}{cc}
      x&y\\
      -yz&-x\\
    \end{array}
  \right),
\end{equation}
where $x^2-y^2z=1$.
\end{lemma}
\begin{pf}
Let $T\in O_{2\times 0+2,\Delta}(\bF)$ and write
$$
T=\left(
    \begin{array}{cc}
      x&y\\
      u&v\\
    \end{array}
  \right),
$$
where $x^2-y^2z=1$, $xu-yvz=0$ and $u^2-v^2z=-z$.
If $xyuv\ne 0$, then $u=x^{-1}yvz$ and $v^2=x^2$, i.e., $v=\pm x$ and $u=\pm yz$. Then $T$ has one of the form (\ref{Eq5}) with $x^2-y^2z=1$.
If $xyuv=0$, then $T$ has one of the following forms
$$
\pm I^{(2)},\,\,
\pm \left(
    \begin{array}{cc}
      1&\\
      &-1\\
    \end{array}
  \right),\,\,
\left(
    \begin{array}{cc}
      0&y\\
      y^{-1}&0\\
    \end{array}
  \right),\,\,
\left(
    \begin{array}{cc}
      0&y\\
      -y^{-1}&0\\
    \end{array}
  \right)
$$
with $-y^2z=1$, which are of the form (\ref{Eq5}). \hfill $\Box$
\end{pf}

Pick a fixed subset $\Omega$ of $\bF^{\ast}$ such that $\bF^{\ast}=\Omega\cup -\Omega$, where $-\Omega=\{-a|a\in\Omega\}$.
Let $E_i$ denote the $\nu$-dimensional column vector having 1 as its $i$-entry and other entries 0's.
Similar to \cite[Theorem~4.1]{Li1}, any element of $\Lambda$ is of the form
$(X-2^{-1}Z\Delta Z^{t}\ \ I^{(\nu)}\ \ Z)$,
where $X\in{\sf K}_{\nu}$ and $Z\in M_{\nu 2}(\bF)$.
Note that $\{\varphi_{0}\}=\{P_{1}\}$ is the trivial orbit of $G_{01}$ on $\Lambda$. We have

\begin{thm}\label{thm2.4}
The nontrivial orbits of $G_{01}$ on $\Lambda$ have the following representatives:
\begin{align}
&\varphi_{1}(r)=([\mathcal{A}_{2r},0^{(\nu-2r)}]\ \ I^{(\nu)}\ \ 0^{(\nu,2)})&(1\leq r\leq \lfloor\nu/2\rfloor),
  \label{E:lwgm1} \\
&\varphi_{2}(r,a)=([-2^{-1}(1-za^2),\mathcal{A}_{2r},0^{(\nu-2r-1)}]\ \ I^{(\nu)}\ \ (E_1\ \ aE_1)) &(0\leq r\leq \lfloor(\nu-1)/2\rfloor),
  \label{E:lwgm2}\\
&\varphi_{3}(r,a)=([\mathcal{A}_{2r},0^{(\nu-2r)}]+[-2^{-1}(1-za^2),0^{(\nu-1)}]\ \ I^{(\nu)}\ \ (E_1\ \ aE_1)) &(1\leq r\leq \lfloor\nu/2\rfloor),
  \label{E:lwgm3}\\
&\varphi_{4}(r)=([0,\mathcal{A}_{2r},0^{(\nu-1-2r)}]+[-2^{-1}\Delta,0^{(\nu-2)}]\ \ I^{(\nu)}\ \ (E_1\ \ E_2)) &(0\leq r\leq \lfloor(\nu-1)/2\rfloor),
  \label{E:lwgm4}\\
&\varphi_{5}(r)=([Y,\mathcal{A}_{2r-2},0^{(\nu-1-2r)}]+[-2^{-1}\Delta,0^{(\nu-2)}]\ \ I^{(\nu)}\ \ (E_1\ \ E_2)) &(1\leq r\leq \lfloor(\nu-1)/2\rfloor),
  \label{E:lwgm5}\\
&\varphi_{6}(r)=([-2^{-1}\Delta,\mathcal{A}_{2r},0^{(\nu-2r-2)}]\ \ I^{(\nu)}\ \ (E_1\ \ E_2)) &(1\leq r\leq \lfloor(\nu-2)/2\rfloor),
  \label{E:lwgm6}\\
&\varphi_{7}(r,b)=([b\mathcal{A}_{2}-2^{-1}\Delta,\mathcal{A}_{2r-2},0^{(\nu-2r)}]\ \ I^{(\nu)}\ \ (E_1\ \ E_2)) &(1\leq r\leq \lfloor\nu/2\rfloor),
  \label{E:lwgm7}\\
&\varphi_{8}(r)=([\mathcal{K},\mathcal{A}_{2r-4},0^{(\nu-2r)}]+[-2^{-1}\Delta,0^{(\nu-2)}]\ \ I^{(\nu)}\ \ (E_1\ \ E_2)) &(2\leq r\leq\lfloor\nu/2\rfloor),
  \label{E:lwgm8}
\end{align}
where $a\in\{0,1\}$, $b\in\Omega$,
$$
Y=\left(
    \begin{array}{ccc}
      0&0&1\\
      0&0&1\\
      -1&-1&0\\
    \end{array}
  \right),
$$
and $\mathcal{K}$ is given by Theorem~\ref{thm2.1}. Moreover the rank of $G_0$ on $\Lambda$ is
$$
(q+7)/2\cdot \lfloor\nu/2\rfloor+4\cdot \lfloor(\nu-1)/2\rfloor+\lfloor(\nu-2)/2\rfloor+3.
$$
\end{thm}
\begin{pf}
Suppose $P\in\Lambda\setminus \{P_1\}$. Then $P=(X-2^{-1}Z\Delta Z^{t}\ \ I^{(\nu)}\ \ Z)$,
where $\mathrm{rank}\,(X-2^{-1}Z\Delta Z^{t}\ \ Z)>0$.

If $Z=0$, then $\mathrm{rank}\,X=2r>0$, which implies that there
exists a $T\in GL_{\nu}(\bF)$ satisfying
$T^{t}XT=[\mathcal{A}_{2r},0^{(\nu-2r)}]$.
Observe $[T,(T^t)^{-1},I^{(2)}]\in G_{01}$ carries $P$ to \eqref{E:lwgm1}.

If $Z\ne 0$, then $\mathrm{rank}\,Z=1$ or 2. We distinguish the following two cases.

{\em Case 1:} $\mathrm{rank}\,Z=1$. Then there exists an $S\in GL_{\nu}(\bF)$
such that $S^{t}Z=(xE_1\ \ yE_1)$, where $(x,y)\ne (0,0)$.
By Theorem~\ref{thm2.1}
there exists a $T\in\mathcal{O}_1$, which is of the form (\ref{Eq2}),
such that $T^{t}(S^{t}XS)T$ is $0^{(\nu)}$ or of form (\ref{Eq3}).
By Lemma~\ref{lem2.2}, there exists an $S_{11}\in O_{2\times 0+2,\Delta}(\bF)$ such that
$T^{t}S^{t}ZS_{11}=b(E_1\ \ aE_1)$, where $a\in\{0,1\}$ and $b\in\bF^{\ast}$.
Observe
$$
(ST)^t(Z\Delta Z^{t})ST=(T^{t}S^{t}ZS_{11})\Delta (S_{11}^{t}Z^{t}ST)
=b^2(E_1\ \ aE_1)\Delta (E_1\ \ aE_1)^t
=[b^2(1-za^2),0^{(\nu-1)}].
$$
Let $T_1=[b^{-1},I^{(\nu-1)}]$. Then
$[STT_1,((STT_1)^t)^{-1},S_{11}]\in G_{01}$ carries $P$ to
$$
((STT_1)^t(X-2^{-1}Z\Delta Z^{t})STT_1\ \ I^{(\nu)}\ \ (STT_{1})^tZS_{11}).
$$
Note that
$\mathrm{rank}((STT_1)^t(X-2^{-1}Z\Delta Z^{t})STT_1\ \ (STT_1)^tZS_{11})=\mathrm{rank}(X-2^{-1}Z\Delta Z^{t}\ \ Z)$.

If $(ST)^{t}X(ST)=0^{(\nu)}$, then $[STT_1,((STT_1)^t)^{-1},S_{11}]$
carries $P$ to \eqref{E:lwgm2} for $r=0$.

If $(ST)^{t}X(ST)=[0,\mathcal{A}_{2r},0^{(\nu-1-2r)}]$, then $[STT_1,((STT_1)^t)^{-1},S_{11}]$
carries $P$ to \eqref{E:lwgm2} for $r>0$.

If $(ST)^{t}X(ST)=[\mathcal{A}_{2r},0^{(\nu-2r)}]$, then $[STT_1,((STT_1)^t)^{-1},S_{11}]$
carries $P$ to
$$
([b^{-1}\mathcal{A}_{2},\mathcal{A}_{2r-2},0^{(\nu-2r)}]+[-2^{-1}(1-za^2),0^{(\nu-1)}]\ \ I^{(\nu)}\ \ (E_1\ \ aE_1)).
$$
Suppose $T_{2}=[1,b,I^{(\nu-2)}]$. Then $[STT_1T_{2},((STT_1T_{2})^t)^{-1},S_{11}]$
carries $P$ to \eqref{E:lwgm3}.

{\em Case 2:} $\mathrm{rank}\,Z=2$.
Then there exists an $S\in GL_{\nu}(\bF)$
such that $S^{t}Z=(E_1\ \ E_2)$.
By Theorem~\ref{thm2.1}, there exists a $T\in\mathcal{O}_2$, which is of the form (\ref{Eq2}) satisfying
$T^{t}(S^{t}XS)T$ is $0^{(\nu)}$ or of form (\ref{Eq4}).
Let $T_1=[T_{11}^{-1},I^{(\nu-2)}]$. Then
$[STT_1,((STT_1)^t)^{-1},I^{(2)}]\in G_{01}$ carries $P$ to
$$
((STT_1)^t(X-2^{-1}Z\Delta Z^{t})STT_1\ \ I^{(\nu)}\ \ (E_1\ \ E_2)).
$$
Observe
$$
(STT_1)^t(Z\Delta Z^{t})STT_1=(T_1^{t}T^{t}S^{t}Z)\Delta (Z^{t}STT_1)=[\Delta,0^{(\nu-2)}]
$$
and
$$
\mathrm{rank}((STT_1)^t(X-2^{-1}Z\Delta Z^{t})STT_1\ \ (E_1\ \ E_2))=\mathrm{rank}(X-2^{-1}Z\Delta Z^{t}\ \ Z).
$$

If $(ST)^{t}X(ST)=0^{(\nu)}$, then $[STT_1,((STT_1)^{t})^{-1},I^{(2)}]$
carries $P$ to \eqref{E:lwgm4} for $r=0$.

If $(ST)^{t}X(ST)=[0,\mathcal{A}_{2r},0^{(\nu-1-2r)}]$, then $[STT_1,((STT_1)^{t})^{-1},I^{(2)}]$
carries $P$ to
$$
([Y_{u,v},\mathcal{A}_{2r-2},0^{(\nu-1-2r)}]+[-2^{-1}\Delta,0^{(\nu-2)}]\ \ I^{(\nu)}\ \ (E_1\ \ E_2)),
$$
where
$$
Y_{u,v}=\left(
    \begin{array}{ccc}
      0&0&u\\
      0&0&v\\
      -u&-v&0\\
    \end{array}
  \right),
\quad
\left(
    \begin{array}{c}
      u\\
      v\\
    \end{array}
  \right)=(T_{11}^{-1})^t\left(
    \begin{array}{c}
      0\\
      1\\
    \end{array}
  \right).
$$
Take $T_{2}=[I^{(2)},v^{-1},I^{(\nu-3)}]$ or $[I^{(2)},u^{-1},I^{(\nu-3)}]$ according to $u=0$ or not, respectively. Then $[STT_1T_{2},((STT_1T_{2})^{t})^{-1},I^{(2)}]$
carries $P$ to \eqref{E:lwgm4} for $r>0$,
or
$$
([Y_{1,c},\mathcal{A}_{2r-2},0^{(\nu-1-2r)}]+[-2^{-1}\Delta,0^{(\nu-2)}]\ \ I^{(\nu)}\ \ (E_1\ \ E_2)),
$$
where $c=u^{-1}v$.
When $c^2-z$ is a square element, we may choose an $s\in\bF^{\ast}$ such that $c^2-z=s^2$.
Let $A=[A_{11},s^{-1},I^{(\nu-3)}]$, where
$$
A_{11}=s^{-1}\left(
    \begin{array}{cc}
      c&-z\\
      -1&c\\
    \end{array}
  \right).
$$
By Lemma~\ref{lem2.3}, $A_{11}^t\in O_{2\times 0+2,\Delta}(\bF)$, and $[A,(A^{t})^{-1},(A_{11}^t)^{-1}]\in G_{01}$ carries
$$
([Y_{1,c},\mathcal{A}_{2r-2},0^{(\nu-2r-1)}]+[-2^{-1}\Delta,0^{(\nu-2)}]\ \ I^{(\nu)}\ \ (E_1\ \ E_2))
$$
to \eqref{E:lwgm4} for $r>0$. When $c^2-z$ is a non-square element, we may choose an $s\in\bF^{\ast}$ such that $s^2(c^2-z)=1-z$.
Let $B=[B_{11},s,I^{(\nu-3)}]$, where
$$
B_{11}=\frac{1}{s(c^2-z)}
\left(
    \begin{array}{cc}
      c-z&z(c-1)\\
      c-1&c-z\\
    \end{array}
  \right).
$$
By Lemma~\ref{lem2.3}, $B_{11}^t\in O_{2\times 0+2,\Delta}(\bF)$, and $[B,(B^{t})^{-1},(B_{11}^t)^{-1}]\in G_{01}$ carries
$$
([Y_{1,c},\mathcal{A}_{2r-2},0^{(\nu-1-2r)}]+[-2^{-1}\Delta,0^{(\nu-2)}]\ \ I^{(\nu)}\ \ (E_1\ \ E_2))
$$
to  \eqref{E:lwgm5}.

If $(ST)^{t}X(ST)=[0^{(2)},\mathcal{A}_{2r},0^{(\nu-2-2r)}]$, then $[STT_1,((STT_1)^{t})^{-1},I^{(2)}]$
carries $P$ to \eqref{E:lwgm6}.

If $(ST)^{t}X(ST)=[\mathcal{A}_{2r},0^{(\nu-2r)}]$, then $[STT_1,((STT_1)^{t})^{-1},I^{(2)}]$
carries $P$ to $\varphi_{7}(r,b)$,
where $(T_{11}^{-1})^{t}\mathcal{A}_2T_{11}^{-1}=b\mathcal{A}_2$ for some $b\in\bF^{\ast}$. Note that $[-1,I^{(\nu-1)},-1,I^{(\nu-1)},-1,1]\in G_{01}$ carries $\varphi_{7}(r,b)$ to  $\varphi_{7}(r,-b)$. So we may choose $b\in\Omega$.

If $(ST)^{t}X(ST)=[\mathcal{K},\mathcal{A}_{2r-4},0^{(\nu-2r)}]$, then
$[STT_1,((STT_1)^{t})^{-1},I^{(2)}]$ carries $P$ to
$$
([U,\mathcal{A}_{2r-4},0^{(\nu-2r)}]+[-2^{-1}\Delta,0^{(\nu-2)}]\ \ I^{(\nu)}\ \ (E_1\ \ E_2)),
$$
where
$$
U=\left(
    \begin{array}{cc}
      0&(T_{11}^{-1})^{t}\\
      -T_{11}^{-1}&0\\
    \end{array}
  \right).
$$
Pick $T_{3}=[I^{(2)},T_{11}^{t},I^{(\nu-4)}]$. Then
$[STT_1T_3,((STT_1T_3)^{t})^{-1},I^{(2)}]$
carries $P$ to  \eqref{E:lwgm8}.

What is left to show that no two subspaces in \eqref{E:lwgm1} -
\eqref{E:lwgm8} can fall into the same orbit.
As an example, we show that any two distinct $\varphi_{2}(r,a)$ and $\varphi_{2}(r,a^{\prime})$ can't fall into the same orbit,
and the rest cases may be handled in a similar way. If there exists an element of $G_{01}$ of form (\ref{Eq1}) carrying
$\varphi_{2}(r,a)$ to $\varphi_{2}(r,a^{\prime})$,
then $T$ is of the form
$$
T=\left(
  \begin{array}{cc}
    t&0\\
    T_{21}&T_{22}\\
  \end{array}
\right),
$$
where $t(1,a)S=(1,a^{\prime})$. By \cite[Theorem~6.4]{Wan}, the subspaces $(1,a)$ and $(1,a^{\prime})$ is of the same type. Since $a,a^{\prime}\in\{0,1\}$, we have $a=a^{\prime}$, a contradiction.

Therefore, the desired result follows. \hfill $\Box$
\end{pf}

For each vertex $Q$ of $\Lambda$, the symbol $\overline Q$ denotes the suborbit containing $Q$.
By \cite[Theorem~1.6, Theorem~3.16, Theorem~6.21]{Wan}, we have
$$
\begin{array}{rcl}
|GL_{\nu}(\bF)|&=&q^{\nu(\nu-1)/2}\prod\limits_{i=1}^{\nu}(q^{i}-1),\\
|Sp_{2\nu}(\bF)|&=&q^{\nu^2}\prod\limits_{i=1}^{\nu}(q^{2i}-1),\\
|O_{2\nu+2,\Delta}(\bF)|&=&q^{\nu(\nu+1)}\prod\limits_{i=1}^{\nu}(q^{i}-1)\prod\limits_{i=0}^{\nu+1}(q^{i}+1).\\
\end{array}
$$

\begin{thm}\label{thm2.5}
The nontrivial orbits of $G_{01}$ on $\Lambda$ have lengths as following:
$$
\begin{array}{rcl}
|\overline{\varphi_{1}(r)}|  &=&\dfrac{|GL_{\nu}(\bF)|}{|Sp_{2r}(\bF)|\cdot |GL_{\nu-2r}(\bF)| \cdot q^{2r(\nu-2r)}},\\[10pt]
|\overline{\varphi_{2}(r,a)}|&=&\dfrac{(q+1)|GL_{\nu}(\bF)|}{|Sp_{2r}(\bF)|\cdot |GL_{\nu-2r-1}(\bF)| \cdot 2q^{(2r+1)(\nu-2r-1)}},\\[10pt]
|\overline{\varphi_{3}(r,a)}|&=&\dfrac{(q+1)|GL_{\nu}(\bF)|}{|Sp_{2r-2}(\bF)|\cdot |GL_{\nu-2r}(\bF)| \cdot 2q^{2r(\nu-2r)+2r-1}},\\[10pt]
|\overline{\varphi_{4}(r)}|  &=&\dfrac{(q+1)\cdot|GL_{\nu}(\bF)|}{|Sp_{2r-2}(\bF)|\cdot |GL_{\nu-2r-1}(\bF)| \cdot 2q^{(2r+1)(\nu-2r)-2}},\\[10pt]
|\overline{\varphi_{5}(r)}|  &=&\dfrac{(q+1)|GL_{\nu}(\bF)|}{|Sp_{2r-2}(\bF)|\cdot |GL_{\nu-2r-1}(\bF)| \cdot 2q^{(2r+1)\nu-4(r^2+1)}},\\[10pt]
|\overline{\varphi_{6}(r)}|  &=&\dfrac{|GL_{\nu}(\bF)|}{|Sp_{2r}(\bF)|\cdot |GL_{\nu-2r-2}(\bF)| \cdot q^{(2r+2)(\nu-2r-2)}},\\[10pt]
|\overline{\varphi_{7}(r,b)}|&=&\dfrac{2\,|GL_{\nu}(\bF)|}{|Sp_{2r-2}(\bF)|\cdot |GL_{\nu-2r}(\bF)| \cdot q^{2r(\nu-2r)}},\\[10pt]
|\overline{\varphi_{8}(r)}|  &=&\dfrac{|GL_{\nu}(\bF)|}{|Sp_{2r-2}(\bF)|\cdot |GL_{\nu-2r}(\bF)| \cdot q^{2r(\nu-2r)+4r-5}}.\\[8pt]
\end{array}
$$
\end{thm}
\begin{pf}
We only calculate
$|\overline{\varphi_{3}(r,a)}|$ and $|\overline{\varphi_{7}(r,b)}|$. The length of other suborbits may be computed in a similar way.

Let $G_{3}(r,a)$ be the stabilizer of $\varphi_{3}(r,a)$ in
$G_{01}$, and let $[T,({T^t})^{-1},S]$ be any element of $G_{3}(r,a)$. Then
$$
T^{t}([\mathcal{A}_{2r},0^{(\nu-2r)}]-[2^{-1}(1-za^2),0^{(\nu-1)}])T=[\mathcal{A}_{2r},0^{(\nu-2r)}]-[2^{-1}(1-za^2),0^{(\nu-1)}]
$$
and $T^{t}(E_1\ \ aE_1)S=(E_1\ \ aE_1)$, which imply that $\mu (1\ \ a)S=(1\ \ a)$ and
$$
 T=\bordermatrix{&_{1}&_{1}&_{2r-2}&_{\nu-2r}\cr
 &\mu&0&0&0\cr
 &t&\mu&-\mu T_{31}^{t}\mathcal{A}_{2r-2}T_{33}&0\cr
 &T_{31}&0&T_{33}&0\cr
 &T_{41}&T_{42}&T_{43}&T_{44}}\hspace{-3pt}
\begin{array}{c}
 _{1}\\
 _{1}\\
 _{2r-2}\\
 _{\nu-2r}
\end{array},
$$
where $\mu^2=1$ and $T_{33}^{t}\mathcal{A}_{2r-2}T_{33}=A_{2r-2}$. By Lemma~\ref{lem2.3}, $\mu (1\ \ a)S=(1\ \ a)$ implies that
$S$ is one of the following forms
$$
\left\{\begin{array}{ll}
\mu I^{(2)},\,\mu\left(
                   \begin{array}{cc}
                     1&0 \\
                     0&-1\\
                   \end{array}
                 \right)& \mathrm{if}\,\,a=0,\\[10pt]
\mu I^{(2)},\,\frac{\mu}{1-z}\left(
                   \begin{array}{cc}
                     1+z&2 \\
                     -2z&-(1+z)\\
                   \end{array}
                 \right)& \mathrm{if}\,\,a=1.\\[10pt]
\end{array}\right.
$$
Hence $|G_{3}(r,a)|=|Sp_{2r-2}(\bF)|\cdot |GL_{\nu-2r}(\bF)|\cdot 4q^{2r(\nu-2r)+2r-1}$ and
$$
|\overline{\varphi_{3}(r,a)}|=[G_{01}:G_{3}(r,a)]
=\frac{(q+1)|GL_{\nu}(\bF)|}{|Sp_{2r-2}(\bF)|\cdot |GL_{\nu-2r}(\bF)| \cdot 2q^{2r(\nu-2r)+2r-1}}.
$$

Let $G_{7}(r,b)$ be the stabilizer of  $\varphi_{7}(r,b)$ in
$G_{01}$. Then $G_{7}(r,b)$ consists of matrices $[T,({T^t})^{-1},(T_{11}^t)^{-1}]$, where
$$
 T=\bordermatrix{&_{2}&_{2r-2}&_{\nu-2r}\cr
 &T_{11}&0&0\cr
 &0&T_{22}&0\cr
 &T_{31}&T_{32}&T_{33}}\hspace{-3pt}
\begin{array}{c}
 _{2}\\
 _{2r-2}\\
 _{\nu-2r}
\end{array},
$$
$T_{11}^t\in O_{2\times 0+2,\Delta}(\bF)$, $T_{11}^t\mathcal{A}_{2}T_{11}=\mathcal{A}_{2}$ and $T_{22}^t\mathcal{A}_{2r-2}T_{2}=\mathcal{A}_{2r-2}$.
By Lemma~\ref{lem2.3}, the matrix $T_{11}^t$ satisfying  $T_{11}^t\in O_{2\times 0+2,\Delta}(\bF)$ and
$T_{11}^t\mathcal{A}_{2}T_{11}=\mathcal{A}_{2}$ is of the form
$$
T_{11}^t=\left(
    \begin{array}{cc}
      x&y\\
      yz&x\\
    \end{array}
  \right),
$$
where $x^2-y^2z=1$.
By \cite[Lemma~1.28]{Wan}, the number of solutions of the equation $x^2-y^2z=1$ is $q+1$. Hence
$|G_{7}(r,b)|=|Sp_{2r-2}(\bF)|\cdot |GL_{\nu-2r}(\bF)|\cdot (q+1)q^{2r(\nu-2r)}$
and
$$
|\overline{\varphi_{7}(r,b)}|=[G_{01}:G_{7}(r,b)]=\frac{2\,|GL_{\nu}(\bF)|}{|Sp_{2r-2}(\bF)|\cdot |GL_{\nu-2r}(\bF)| \cdot q^{2r(\nu-2r)}}.
$$
\hfill $\Box$
\end{pf}

\section{Quasi-strongly
regular graphs}

As a generalization of strongly regular graphs, quasi-strongly
regular graphs were discussed by W. Golightly, W. Haynworth and D.G. Sarvate \cite{gol1} and
F. Goldberg \cite{Gold}.
Let $c_1,c_2,\dots,c_p$ be distinct non-negative integers. A connected graph of degree $k$ on $n$ vertices is
{\em quasi-strongly regular} with parameters
$(n,k,\lambda;c_1,c_2,\dots,c_p)$ if   any two adjacent vertices
have $\lambda$ common neighbors, and any two non-adjacent vertices
have $c_i$ common neighbors  for some $i$ $(1\le i \le p)$.

Since $\Gamma$ is a regular near polygon,
the induced subgraph on $\Lambda$ is edge regular, denoted by the same symbol $\Lambda$. Therefore, $\Lambda$ is
quasi-strongly regular. In this section we compute all the
parameters of $\Lambda$.

Let $\Lambda(P)$ be the set of neighbors of $P$ in $\Lambda$.
Clearly, $|\Lambda(P)\cap\Lambda(Q)|=0$ whenever $\dim\,(P+Q)>\nu+2$.
Note that for the vertex $P_1$ in $\Lambda$ as in Section 2, the subspace $Q\in\Lambda$ satisfying $\dim\,(Q+P_1)=\nu+2$ lies in the set
$$
\overline{\varphi_{1}(1)}\,\cup\,\overline{\varphi_{4}(0)}\,\bigcup\limits_{a\in\{0,1\}}\overline{\varphi_{3}(1,a)}\,\bigcup\limits_{b\in\Omega}\overline{\varphi_{7}(1,b)}.
$$
To study $|\Lambda(P)\cap\Lambda(Q)|$ for any two vertices $P$
and $Q$ with $\dim\,(P+Q)=\nu+2$, by Theorem~\ref{thm2.4}, it suffices to
consider $|\Lambda(P_1)\cap\Lambda(Q)|$, where $Q\in\{\varphi_{1}(1),\varphi_{3}(1,a),\varphi_{4}(0),\varphi_{7}(1,b)\}$, $a\in\{0,1\}$ and $b\in\Omega$.

\begin{lemma}~\label{lem3.1}
For any vertex $R=(A-2^{-1}C\Delta C^t\ \ I^{(\nu)}\ \ C)$ of $\Lambda$, the neighborhood of $R$ is
$$
\Lambda(R)=\{(A-2^{-1}(C\Delta C^t+D\Delta D^t+2D\Delta C^t)\ \ I^{(\nu)}\ \ C+D)\,|\,D \in M_{\nu 2}(\bF),\,\,\mathrm{rank}\,D=1\}.
$$
\end{lemma}
\begin{pf}
Note that $\Lambda(R)$ consists of matrices with the form
$(A-2^{-1}C\Delta C^t+X\ \ I^{(\nu)}\ \ C+D)$,
where $X\in M_{\nu}(\mathbb{F}_{q})$, $D\in M_{\nu 2}(\mathbb{F}_{q})$, $\mathrm{rank}\,(X\ \ D)=1$ and $X+X^t+C\Delta D^t+D\Delta C^t+D\Delta D^t=0$.
It follows that $\mathrm{rank}\,D=1$. So we may write $D=D_0(x\,\,y)$ and
$X=D_0W^t$, where $0\not=D_0\in M_{\nu 1}(\mathbb{F}_q)$, $(x,\,y)\not=(0,\,0)$ and
$W\in M_{\nu 1}(\mathbb{F}_q)$.
Let $C=(C_1\,\,C_2)$ and $TD_0=E_{1}$ for some $T\in GL_{\nu}(\mathbb{F}_q)$. Then
$$
E_{1}(T(W\!+\!xC_1\!-\!yzC_2))^t\!+\!T(W\!+\!xC_1\!-\!yzC_2)E_{1}^t\!+\!(x^2-y^2z)E_{1}E_{1}^t=0.
$$
It follows that $T(W+xC_1-yzC_2)=-2^{-1}(x^2-y^2z)E_{1}$. So $W=-2^{-1}(x^2-y^2z)D_0-xC_1+yzC_2$ and
$X=-2^{-1}(D{\Delta}D^t+2D\Delta C^t)$. The desired result follows.  \hfill $\Box$
\end{pf}

Note that when $\nu=1$, any element of $\Lambda$ is of the form
$(-2^{-1}(a^2-zb^2)\ \ 1\ \ (a\ \ b))$,
where $a,b\in\bF$. Then $\Lambda$ is a clique with $q^2$ vertices.

\begin{lemma}~\label{lem3.2}
Let $P$ and $Q$ be any two vertices of $\Lambda$ with $\dim\,(P+Q)=\nu+2$. If $\nu\geq 2,$
then $|\Lambda(P)\cap\Lambda(Q)|$ is equal to
$0$, $q^2$, $q^2-1$ or $q^2+q$.
\end{lemma}
\begin{pf}
For any $Q\in\{\varphi_{1}(1), \varphi_{3}(1,a),\varphi_{4}(0),\varphi_{7}(1,b)\}$, it suffices to show that $|\Lambda(P_1)\cap\Lambda(Q)|=0,q^2,q^2-1$ or $q^2+q$. We only compute $|\Lambda(P_1)\cap\Lambda(\varphi_{3}(1,a))|$, and the others can be treated similarly.

Let $R\in\Lambda(P_1)\cap\Lambda(\varphi_{3}(1,a))$. From $R\in\Lambda(P_1)$ and Lemma~\ref{lem3.1} we know that $R$ is of the form $R=(-2^{-1} D\Delta D^t\ \ I^{(\nu)}\ \ D)$, where $D\in M_{\nu 2}(\bF)$ and $\mathrm{rank}\,D=1$. Again
from $R\in\Lambda(\varphi_{3}(1,a))$ and Lemma~\ref{lem3.1} we know that
$$
R=([\mathcal{A}_{2},0^{(\nu-2)}]-[2^{-1}(1-za^2),0^{(\nu-1)}]-2^{-1}(D_1\Delta D_1^t+2D_1\Delta (E_1\ \ aE_1)^t)\ \ I^{(\nu)}\ \ (E_1\ \ aE_1)+D_1),
$$
where $D_1\in M_{\nu 2}(\bF)$ and $\mathrm{rank}\,D_1=1$. Therefore, $D=(E_1\ \ aE_1)+D_1$ and
$$
-2^{-1}D\Delta D^t=[\mathcal{A}_{2},0^{(\nu-2)}]-[2^{-1}(1-za^2),0^{(\nu-1)}]-2^{-1}(D_1\Delta D_1^t+2D_1\Delta (E_1\ \ aE_1)^t).
$$
It follows that $-2^{-1}(E_1\ \ aE_1)\Delta D_1^t=[\mathcal{A}_{2},0^{(\nu-2)}]-2^{-1}D_1\Delta (E_1\ \ aE_1)^t$; and so $D$ is of the form
$$
D=\left(
        \begin{array}{ccccc}
          c+1&zad_{2}-2&0&\cdots &0\\
          a+d_{1}&d_{2}&0&\cdots &0\\
        \end{array}
      \right)^t,
$$
where $cd_{2}-d_{1}(zad_{2}-2)=0$.
Observe the number of solutions $(c,d_{1},d_{2})$ satisfying $cd_{2}-d_{1}(zad_{2}-2)=0$ is $q+(q-1)q=q^2$. Hence $|\Lambda(P_1)\cap\Lambda(\varphi_{3}(1,a))|=q^2$.
\hfill $\Box$
\end{pf}

\begin{thm}\label{thm3.3} Let $\nu\geq 2$. Then
$\Lambda$ is a quasi-strongly regular graph with parameters
$$
(q^{\nu(\nu+3)/2},(q^{\nu}-1)(q+1),q^{\nu}+q^2-q-1;0,q^2,q^2-1,q^2+q).
$$
\end{thm}
\begin{pf}
Since $\Lambda$ consists of the vertices as  the form $(X-2^{-1}Z\Delta Z^t\ \ I^{(\nu)}\ \ Z)$, where $X$ is a $\nu\times\nu$ alternate matrix, and $Z\in M_{\nu 2}(\bF)$, we have $n=q^{\nu(\nu+3)/2}$. By Theorem~\ref{thm2.5},
$$
k =|\overline{\varphi_{2}(0,0)}|+|\overline{\varphi_{2}(0,1)}|=2|\overline{\varphi_{2}(0,0)}|=(q^{\nu}-1)(q+1).
$$
Note that $\varphi_{2}(0,a)\in\Lambda(P_1)$.
In order to compute the parameter $\lambda$,
it  suffices to compute the size of the common neighbors of $P_1$ and $\varphi_{2}(0,a)$.
Let $R\in\Lambda(P_1)\cap\Lambda(\varphi_{2}(0,a))$. From $R\in\Lambda(P_1)$ and Lemma~\ref{lem3.1} we know that $R$ is of the form
$R=(-2^{-1}D\Delta D^t\ \ I^{(\nu)}\ \ D)$, where $D\in M_{\nu 2}(\bF)$ and $\mathrm{rank}\,D=1$.
Similar to the proof of Case 2 in Lemma~\ref{lem3.2}, $D$ is of rank 1 with the form
$$
D=\left(
        \begin{array}{cccc}
          c+1&azd_{2}&\cdots &azd_{\nu}\\
          d_{1}+1&d_{2}&\cdots &d_{\nu}\\
        \end{array}
      \right)^t.
$$
Observe the number of matrices $D$ satisfying $(d_{2},\ldots,d_{\nu})=(0,\ldots,0)$ and $(d_{2},\ldots,d_{\nu})\ne (0,\ldots,0)$ are $q^2-1$ and $(q^{\nu-1}-1)q$, respectively. So $\lambda=|\Lambda(P_1)\cap\Lambda(\varphi_{2}(0,a))|=q^{\nu}+q^2-q-1$. The rest parameters of $\Lambda$ are listed in Lemma~\ref{lem3.2}.
\hfill $\Box$
\end{pf}

\section{Association schemes}

In this section we discuss the association scheme based on $\Lambda$
when $\nu=2$.

A {\em $d$-class association scheme} $\mathfrak{X}$ is a pair $(X,\{R_i\}_{0\leq i\leq d})$, where $X$ is a finite set, and each $R_{i}$ is a nonempty subset of $X\times X$ satisfying the following axioms:
\begin{enumerate}[(i)]
\item \label{as1}
$R_0=\{(x,x)\,|\,x\in X\}$;

\item
$X\times X=R_0\cup R_1\cup\cdots\cup R_d$, $R_i\cap
R_j=\emptyset\ (i\not=j)$;

\item
${}^tR_i=R_{i'}$ for some $i'\in\{0,1,\ldots,d\}$, where
${}^tR_i=\{(y,x)\,|\,(x,y)\in R_i\}$;

\item \label{as4}
for all $i,j,k\in\{0,1,\ldots, d\}$, there exists an integer
$p_{ij}^{k}=|\{z\in X\,|\,(x,z)\in R_i,\,(z,y)\in R_j\}|$ for every
$(x,y)\in R_{k}$.
\end{enumerate}
The integers $p_{ij}^{k}$ are called the {\em intersection numbers}
of $\mathfrak{X}$, and $k_{i}\,(=p_{ii^{\prime}}^{0})$ is called the
{\em valency} of $R_i$.
Furthermore, $\mathfrak{X}$ is called {\em symmetric} if $i^{\prime}=i$ for all
$i$.
As for more information concerning association schemes, the readers may
consult \cite{BI,BCN}.

Let $G$ be a transitive permutation group on a finite set $X$, and
$R_{0},R_{1},\ldots,R_{d}$ be the orbits of the induced action of
$G$ on $X\times X$. It is well known that  $(X,\{R_i\}_{0\leq i\leq
d})$ is an association scheme (\cite[\S 2.2]{BI}).

Note that the action of $G_0$ on $\Lambda\times \Lambda$ determines
an association scheme. We shall discuss the association scheme in
the case $\nu=2$.

In the rest we always assume that $\nu=2$. By Theorem~\ref{thm2.4},
the orbits of $G_{01}$ on $\Lambda$ have the following
representatives:
$$
\varphi_{0},\varphi_{1}(1),\varphi_{2}(0,a),\varphi_{3}(1,a),\varphi_{4}(0),\varphi_{7}(1,b),
$$
where $a\in\{0,1\}$, $b\in\Omega$. For the action of $G_{0}$ on
$\Lambda\times\Lambda$, let
$R_{0},R_{1},R_{2_{a}},R_{3_{a}},R_{4},R_{5_{b}}$ denote the orbits
containing $(\varphi_{0},\varphi_{0})$,
$(\varphi_{0},\varphi_{1}(1))$, $(\varphi_{0},\varphi_{2}(0,a))$,
$(\varphi_{0},\varphi_{3}(1,a))$, $(\varphi_{0},\varphi_{4}(0))$,
$(\varphi_{0},\varphi_{7}(1,b))$, respectively. Then
$R_{0},R_{1},R_{2_{a}},R_{3_{a}},R_{4},R_{5_{b}}$ are all the orbits
of the action of $G_{0}$ on $\Lambda\times\Lambda$.

Let $G_{\varphi_{1}(1)}$ be the stabilizer of $\varphi_{1}(1)$ in
$G_{01}$. Then $G_{\varphi_{1}(1)}$ consists of matrices with the
form $[T,(T^{t})^{-1},S]$, where $T^{t}\in Sp_{2}(\mathbb{F}_q)$ and
$S\in O_{2\times 0+2,\Delta}(\mathbb{F}_q)$. So
$$
|G_{\varphi_{1}(1)}|=|Sp_{2}(\mathbb{F}_{q})|\cdot|O_{2\times 0+2,\Delta}(\mathbb{F}_{q})|=2q(q-1)(q+1)^2.
$$

In order to discuss the association scheme based on $\Lambda$, we need the following lemmas.

\begin{lemma}~\label{lem4.1}
The orbits of $G_{\varphi_{1}(1)}$ on $\Lambda$ have the following representatives:
$$
\begin{array}{rcl}
\phi_{1x}&=&(x\mathcal{A}_{2}\ \ I^{(2)}\ \ 0^{(2)}),\\[0.1cm]
\phi_{2x,a}&=&\big(x\mathcal{A}_{2}+[-\frac{1}{2}(1-za^{2}),0]\ \ I^{(2)}\ \ (E_{1}\ \ aE_{1})\big),\\[0.1cm]
\phi_{3x,c}&=&\big(x\mathcal{A}_{2}+[-\frac{1}{2}c^{2},\frac{1}{2}z]\ \ I^{(2)}\ \ [c,1]\big),\\
\end{array}
$$
where $a\in\{0,1\}$, $x\in\mathbb{F}_{q}$ and $c\in\Omega$.
\end{lemma}
\begin{pf}
The proof is similar to that of Theorem~\ref{thm2.4}, and is omitted.
\hfill $\Box$
\end{pf}

\begin{lemma}~\label{lem4.2}
Let $c\in\Omega$. Then the number of $(T,S)$ satisfying $T^{t}\in Sp_{2}(\mathbb{F}_{q})$, $S\in O_{2\times 0+2,\Delta}(\mathbb{F}_q)$ and $T^{t}[c,1]S=[c,1]$ is $q+1$.
\end{lemma}
\begin{pf}
Since $T^{t}[c,1]S=[c,1]$, by Lemma~\ref{lem2.3},
$$
T=\left(\begin{array}{cc}
\mu&c^{-2}sz\\
s&\mu\\
\end{array}\right),
S=\left(
    \begin{array}{cc}
      \mu&-c^{-1}s\\
      -c^{-1}sz&\mu\\
    \end{array}
  \right),
$$
where $\mu,s,c\in\mathbb{F}_{q}$ and $\mu^2-c^{-2}s^2z=1$. By \cite[Lemma~1.28]{Wan}, the number of $(\mu,c)$ satisfying $\mu^2-c^{-2}s^2z=1$ is $q+1$, as desired.
\hfill $\Box$
\end{pf}

\begin{lemma}~\label{lem4.3}
The representatives $\phi_{1x},\phi_{2x,a},\phi_{3x,c}$ listed in Lemma~\ref{lem4.1} satisfy
\begin{align*}
\begin{array}{ll}
(\varphi_{0},\phi_{10})\in R_{0},&(\phi_{10},\varphi_{1}(1))\in R_{1},\\
(\varphi_{0},\phi_{11})\in R_{1},&(\phi_{11},\varphi_{1}(1))\in R_{0},\\
(\varphi_{0},\phi_{1d})\in R_{1},&(\phi_{1d},\varphi_{1}(1))\in R_{1},\\[0.1cm]
(\varphi_{0},\phi_{20,a})\in R_{2_{a}},&(\phi_{20,a},\varphi_{1}(1))\in R_{3_{a}},\\
(\varphi_{0},\phi_{21,a})\in R_{3_{a}},&(\phi_{21,a},\varphi_{1}(1))\in R_{2_{a}},\\
(\varphi_{0},\phi_{2d,a})\in R_{3_{a}},&(\phi_{2d,a},\varphi_{1}(1))\in R_{3_{a}},\\[0.1cm]
(\varphi_{0},\phi_{30,c})\in R_{4},&(\phi_{30,c},\varphi_{1}(1))\in R_{5_{\varepsilon c^{-1}}},\\
(\varphi_{0},\phi_{31,c})\in R_{5_{\varepsilon c^{-1}}},&(\phi_{31,c},\varphi_{1}(1))\in R_{4},\\
(\varphi_{0},\phi_{3d,c})\in R_{5_{\varepsilon_{1} c^{-1}d}},&(\phi_{3d,c},\varphi_{1}(1))\in R_{5_{\varepsilon_{2} c^{-1}(1-d)}},\\
\end{array}
\end{align*}
where $d\in\mathbb{F}_{q}\backslash\{0,1\}$,
$\varepsilon,\varepsilon_{1},\varepsilon_{2}\in\{1,-1\}$, $\varepsilon c^{-1},\varepsilon_{1} c^{-1}d,\varepsilon_{2} c^{-1}(1-d)\in\Omega$.
\end{lemma}
\begin{pf}
We only show $(\varphi_{0},\phi_{30,c})\in R_{4}$ and $(\phi_{30,c},\varphi_{1}(1))\in R_{5_{\varepsilon c^{-1}}}$. The left cases may be treated similarly, and will be omitted.
Note that
$[c^{-1},1,c,1,I^{(2)}]\in G_{0}$
carries $\varphi_{0}$ and $\phi_{30,c}$ to $\varphi_{0}$ and $\varphi_{4}(0)$,
respectively, so $(\varphi_{0},\phi_{30,c})\in R_4$. Let $\varepsilon=1$ or $-1$ according to $c^{-1}\in\Omega$ or $-c^{-1}\in\Omega$, respectively.
Then
$$
\left(\begin{array}{cccccc}
                -c^{-1}&0&&&&\\
                0&-\varepsilon&&&&\\
                \frac{1}{2}c&0&-c&0&-c&0\\
                0&-\frac{1}{2}\varepsilon z&0&-\varepsilon&0&-\varepsilon\\
                -1&0&&&1&0\\
                0&\varepsilon z&&&0&\varepsilon\\
                 \end{array}\right)\in G_{0}
$$
carries $\phi_{30,c}$ and $\varphi_{1}(1)$ to $\varphi_{0}$ and $\varphi_{7}(1,\varepsilon c^{-1})$, respectively, which implies  $(\phi_{30,c},\varphi_{1}(1))\in R_{5_{\varepsilon c^{-1}}}$.
\hfill $\Box$
\end{pf}

\begin{thm}\label{thm4.4} The configuration $\mathfrak{X}=(\Lambda,\{R_{0},R_{1},R_{2_{a}},R_{3_{a}},R_{4},R_{5_{b}}\}_{a\in\{0,1\},b\in\Omega})$ is a symmetric association scheme with class $(q+11)/2$, whose non-zero intersection numbers $p^{1}_{ij}$ are given by
$$
\begin{array}{ll}
&p^{1}_{01}=p^{1}_{10}=1,\;\;
 p^{1}_{11}=q-2,\;\;
 p^{1}_{2_{a},3_{a}}=p^{1}_{3_{a},2_{a}}=(q-1)(q+1)^2/2,\\[0.1cm]
&p^{1}_{3_{a},3_{a}}=(q-2)(q-1)(q+1)^{2}/2,\;\;
 p^{1}_{4,5_{b}}=p^{1}_{5_{b},4}=p^{1}_{5_{b},5_{\varepsilon bd^{-1}(1-d)}}=2q(q^2-1), \\[0.1cm]
\end{array}
$$
where $d\in\mathbb{F}_{q}\backslash\{0,1\}$, $\varepsilon\in\{1,-1\}$ and $\varepsilon bd^{-1}(1-d)\in\Omega$.
\end{thm}
\begin{pf}
By Theorem~\ref{thm2.4}, $\mathfrak{X}$ forms an association scheme of class $(q+11)/2$.

Now we prove   $\mathfrak{X}$ is symmetric.  Since
$$
\left(
  \begin{array}{ccccc}
    1&0&&&\\
    0&-1&&&\\
    0&1&1&0&\\
    1&0&0&-1&\\
    &&&&I^{(2)}\\
  \end{array}
\right)\in G_{0}
$$
interchanges  $\varphi_{1}(1)$ and $\varphi_{0}$, ${}^t R_1=R_1$.
The left cases can be treated similarly, and will be omitted .

In order to compute non-zero intersection numbers $p^{1}_{ij}$ of $\mathfrak{X}$, we need consider the cases listed in Lemma~\ref{lem4.3}. Here we only calculate $p^{1}_{2_{a},3_{a}}$ and $p^{1}_{4,5_{b}}$ by the way of examples.

Let $G_{\phi_{20,a}}$  be the stabilizer of $\phi_{20,a}$ in $G_{\varphi_{1}(1)}$. By Lemma~\ref{lem4.3}, $p^{1}_{2_{a},3_{a}}=[G_{\varphi_{1}(1)}:G_{\phi_{20,a}}]$, the index of $G_{\phi_{20,a}}$ in $G_{\varphi_{1}(1)}$.
Note that $G_{\phi_{20,a}}$ consists of matrices $[T,(T^{t})^{-1},S]$, where
$$
T^{t}\in Sp_{2}(\mathbb{F}_q),\; T^{t}[1,0]T=[1,0],\; S\in O_{2\times 0+2,\Delta}(\mathbb{F}_q)\;\;\mbox{and}\;\;
T^{t}(E_{1}\ \ aE_{1})S=(E_{1}\ \ aE_{1}).
$$
It follows that
$$
T=\left(
    \begin{array}{cc}
      \mu &0\\
      t&\mu\\
    \end{array}
  \right),\;\;
S=\mu I^{(2)}\;\;\mbox{or}\;\;
S=\frac{\mu}{1-za^{2}}\left(
    \begin{array}{cc}
      1+za^2 &2a\\
      -2za&-(1+za^2)\\
    \end{array}
  \right),
$$
where  $\mu^{2}=1$. Therefore, $|G_{\phi_{20,a}}|=4q$ and
$$
p^{1}_{2_{a},3_{a}}=[G_{\varphi_{1}(1)}:G_{\phi_{20,a}}]=|G_{\varphi_{1}(1)}|/|G_{\phi_{20,a}}|
=(q-1)(q+1)^2/2.
$$

Let $G_{\phi_{30,c}}$ be the stabilizer of $\phi_{30,c}$ in $G_{\varphi_{1}(1)}$, where $c\in\Omega$.
Then $G_{\phi_{30,c}}$ consists of matrices $[T,(T^{t})^{-1},S]$, where
$$
T^{t}\in Sp_{2}(\mathbb{F}_q),\;
T^{t}[c^2,-z]T=[c^2,-z],\;
S\in O_{2\times 0+2,\Delta}(\mathbb{F}_q)\;\;\mbox{and}\;\; T^{t}[c,1]S=[c,1].
$$
Note that $T^{t}[c,1]S=[c,1]$ implies $T^{t}[c^2,-z]T=[c^2,-z]$. By Lemma~\ref{lem4.2}, $|G_{\phi_{30,c}}|=q+1$; and by Lemma~\ref{lem4.3},
$$
p^{1}_{4,5_{\varepsilon c^{-1}}}=[G_{\varphi_{1}(1)}:G_{\phi_{30,c}}]=|G_{\varphi_{1}(1)}|/|G_{\phi_{30,c}}|
=2q(q^2-1),
$$
which is independent of choices of $c\in\Omega$. Hence $p^{1}_{4,5_{b}}=2q(q^2-1)$.
\hfill $\Box$
\end{pf}

\paragraph*{Remarks}
All the valencies of $\mathfrak{X}$ are given by Theorem~\ref{thm2.5}. By
a similar method in this section, all the
intersection numbers $p_{ij}^{k}$ of $\mathfrak{X}$ can be calculated.

\section*{Acknowledgment}

This research is partially supported by Hunan Provincial Natural
Science Foundation of China(09JJ3006),  NCET-08-0052, NSF of China
(10871027, 10971052),  the Fundamental Research Funds for the
Central Universities of China, and Langfang Teachers' College
(LSZB201005).

\end{document}